\documentclass{amsart}
\usepackage[dvipdfmx]{graphicx}

\usepackage{amssymb}
\usepackage{amsmath}
\usepackage{url}
\usepackage{enumerate}

\usepackage{hyperref}
\hypersetup{
  colorlinks  =   true,
  linkcolor   =   blue,
  filecolor   =   blue,
  urlcolor    =   blue,
  citecolor   =   blue,
}

\begin{document}

\title[Kummer surfaces]{Kummer surfaces:200 years of study}
\author{Igor Dolgachev}

\address{Department of Mathematics,
University of Michigan,
525 East University Avenue,
Ann Arbor, MI 48109-1109 USA}
\email{idolga@umich.edu}


\begin{abstract}
This is a brief history of discovery and later study of Kummer surfaces. The article is based on the author's  Oliver Club talk at Cornell University on October 10, 2019 delivered exactly 101 years since the first talk at the club given by John Hutchinson.
\end{abstract}

\maketitle

The fascinating story about the Kummer surface starts from discovery by Augustin-Jean Fresnel in 1822 of the equation describing the propagation of light in an optically biaxial crystal \cite{Fresnel}.\marginpar{Augustin-Jean Fresnel:1822} Let $v$ denote the light speed in some matter and $n = c/v$, where $c$ is the light speed in vacuum. It depends on the coordinates $x = (x_1,x_2,x_3)$  of a point  and the unit direction vector $\xi = (\xi_1,\xi_2,\xi_3)$. The propagation of light is described by the function $n(x,\xi)$. The matter is homogeneous if $n$ does not depend on $x$ and isotropic if it does not depend on the direction. For example, Maxwell studied the  fisheye and found that the light is inhomogeneous but isotropic. A biaxial crystal gives an example of homogeneous but anisotropic propagation. Fresnel found the equation of propagation of light in such a crystal of the form:
\begin{equation}\label{eq1}
\frac{\xi_1^2}{\frac{1}{a_1^{2}}-\frac{1}{n^{2}}}+\frac{\xi_2^2}{\frac{1}{a_2^{2}}-\frac{1}{n^2}}+\frac{\xi_3^2}{\frac{1}{a_3^{2}}-\frac{1}{n^2}} = 0,
\end{equation}
where $a_1,a_2,a_3$ are constants describing the property of the crystal (principal refraction indices). In new coordinates 
$(x,y,z) = n(\xi_1,\xi_2,\xi_3)$ one can rewrite the equation \eqref{eq1} in homogeneous form as
\begin{equation}
\label{eq2}
\frac{a_1^2x^2}{x^2+y^2+z^2-a_1^2w^2}+\frac{a_2y^2}{x^2+y^2+z^2-a_2^2w^2}+\frac{a_3z^2}{x^2+y^2+z^2-a_3^2w^2} = 0
\end{equation}
After cleating the  denominators, we find a homogeneous equation of degree $4$. 

In 1833 Sir William Hamilton  \marginpar{William\,Hamilton:1833} discovered that the surface has four real singular points
$$(\pm a_3\sqrt{\frac{a_1^2-a_2^2}{a_1^2-a_3^2}},0,\pm a_1\sqrt{\frac{a_2^2-a_3^2}{a_1^2-a_3^2}},1)$$
and also $4$ real \emph{trope-conics} cut out by the planes $\alpha x_1+\beta x_2+\gamma x_3+x_4 = 0$, where 
$$(\alpha,\beta,\gamma,1) = (\pm \frac{a_3}{a_2^2}\sqrt{\frac{a_1^2-a_2^2}{a_1^2-a_3^2}},0,\pm \frac{a_1}{a_2^2}\sqrt{\frac{a_2^2-a_3^2}{a_1^2-a_3^2}},1)$$
\cite[p. 134]{Hamilton}. In fact, over $\mathbb{C}$, it has additionally 12 nodes and 12 trope-conics such that the incidence configuration of $16$ nodes and $16$ trope-conics is of type $(16_6)$ (this the famous \emph{Kummer abstract configuration}).
In 1849 Arthur Cayley \marginpar{Arthur Cayley: 1848} proved that the Fresnel's wave surface is a particular case of a quartic surface. It contains four pairs of coplanar conics,  the three vertices  of  the tetrahedron formed by these four planes lying in the same face  are conjugate with respect to the both conics \cite{Cayley}. He gave the name  \emph{tetraedroid quartic surface} to such a quartic surface.   He also discovered an important property of tetraedroid quartic surfaces: they are projectively \emph{self-dual} (or \emph{reciprocal}). The wave surfaces were the subject of study for many  mathematicians of the 19th century. Among them were  A. Cauchy, A. Cayley, J. Darboux, J. Mac Cullagh, J. Sylvester and W. Hamilton (see \cite[pages 114-115]{Loria}). A nice modern exposition of the theory of Fresnel's wave surfaces can be found in \cite{Knorrer}.

Projective equivalence classes of Fresnel wave surfaces depend on $2$ parameters and as we shall see momentary Kummer surfaces depend on $3$ parameters. 

In 1847 Adolph G\"opel \marginpar{Adolph G\"opel:1847} using the transcendental theory of theta functions had found a relation of order four between theta functions of second order in two variables that expresses an equation of a general Kummer surface \cite{Gopel}. 

Let $T = \mathbb{C}^g/\Lambda$ be a compact $g$-dimensional torus, the quotient of $\mathbb{C}^g$ by the group of translations $\Lambda$ isomorphic to $\mathbb{Z}^{2g}$. There are no holomorphic functions on $T$ because it is a compact complex manifold, instead one considers nonzero holomorphic sections of a holomorphic line bundle $L$ on $T$. For general $\Lambda$ no line bundle has enough sections to embed $T$   into a projective space, i.e. $T$ does not admit any ample (or positive) line bundle. However,  in 1857 Bernard Riemann  found a condition on $\Lambda$ such that such $L$ exists and  $\dim \Gamma(T,L^n) = n^g$. Such complex tori are now called  \emph{principally polarized abelian varieties}.  They depend on $\frac{1}{2}g(g+1)$ parameters. Holomorphic sections of $L^n$ can be lifted to holomorphic functions on $\mathbb{C}^g$ that are $\Lambda$-invariant up to some multiplicative factor. They are called \emph{theta-functions} of order $n$. For $n = 1$, such a holomorphic function  is the famous Riemann theta function $\Theta(z,\Lambda)$. One can modify Riemann's  expression for $\Theta(z,\Lambda)$  to obtain $n^{2g}$  theta functions $\Theta_{m,m'}(z,\Lambda)$ with characteristics $(m,m')\in (\mathbb{Z}/n\mathbb{Z})^g\oplus (\mathbb{Z}/n\mathbb{Z})^g$. They generate  the linear space $\Gamma(T,L^n)$. 

An example of a principally polarized abelian variety is the \emph{Jacobian variety} $\textrm{Jac}(C)$ of a Riemann surface $C$ of genus $g$.  Here $\Lambda$ is spanned by vectors 
$v_i  = (\int_{\gamma_i}\omega_1,\ldots,\int_{\gamma_{i}}\omega_g)\in \mathbb{C}^g,\  i = 1,\ldots,2g,$
where $(\omega_1,\ldots, \omega_g)$ is a basis of the linear space of holomorphic differential $1$-forms and $(\gamma_1,\ldots,\gamma_{2g})$ is a basis of $H_1(T,\mathbb{Z})$.

In the case $g = 2$ and  $n = 2$ ,   G\"opel was able to find a special basis $(\theta_0,\theta_1,\theta_2,\theta_3)$ in the space $\Gamma (\textrm{Jac}(C),L^2) \cong \mathbb{C}^4$ such that the map
$$\Phi:T\to \mathbb{P}^3, (z_1,z_2) \mapsto (\theta_0(z),\theta_1(z),\theta_2(z),\theta_3(z))$$
satisfies $\Phi(-z_1,-z_2) = \Phi(z_1,z_2)$ and its image $X$ is the set of zeros in $\mathbb{P}^3$ of a quartic polynomial
{\small $$F = A(x^4+y^4+z^4+w^4)+2B(x^2y^2+z^2w^2)+2C(x^2z^2+y^2w^2)+2D(x^2w^2+y^2z^2)+4Exyzw,$$}
where the coefficients $(A,B,D,E)$ satisfy a certain explicit equation in terms of theta constants, the values of theta functions at $0$.  The images of sixteen $2$-torsion points $\epsilon\in \frac{1}{2}\Lambda/\Lambda\in \mathbb{C}^2/\Lambda$ are singular points of  $X$. The Abel-Jacobi map $C\to \textrm{Jac}(C), x\mapsto (\int_{x_0}^x\omega_1,\int_{x_0}^x\omega_2)\mod \Lambda$ embeds $C$ into $\textrm{Jac}(C)$ and the images of the curves $C+\epsilon$ are the 16 trope-conics of $X$.

 In 1864 Ernst Kummer \marginpar{Ernst Kummer:1864} had shown that the Fresnel's wave surface represents a special case of a 3-parametrical family of $16$-nodal quartic surfaces \cite{Kummer1}. He proved that they must contain 16 trope-conics which together with $16$ nodes form an abstract incidence configuration $(16_6)$. Kummer shows that any 16-nodal quartic surface has a tetrahedron with conic-tropes in the faces intersecting at two points on the edges. No vertex is a node. From this he deduced that  there exists a quadric surface that contains the four trope-conics. Using this he finds an  equation of a general Kummer surface of the form
$$(x^2+y^2+z^2+w^2+a(xy+zw)+b(xz+yw)+c(xw+yz))^2 + Kxyzw = 0,$$
where $xyzw=0$ is the  equation of a chosen tetrahedron and $K = a^2+b^2+c^2-2abc-1$. Since then any 16-nodal quartic surface became known as a \emph{Kummer quartic surface}. 

It took almost 30 years since G\"opel's discovery to realize that the G\"opel equation, after a linear change of variables, can be reduced to Kummer equation. This was done by Carl Borchardt in 1878 \cite{Borchardt}. \marginpar{Carl Borchardt:1878}

In fact, G\"opel's discovery leads to a modern definition of the Kummer surface and its  higher-dimensional version, the \emph{Kummer variety}. One considers any $g$-dimensional complex torus $T$ and divides it by the involution $z\mapsto -z$, the orbit space is called the \emph{Kummer variety} of $T$ and is denoted by $\textrm{Kum}(T)$. If $T$ admits an embedding into a projective space, then $\textrm{Kum}(T)$ is a projective algebraic variety with $2^{2\dim T}$ singular points. Moreover, if $T$ is a principally polarized abelian variety, for example the Jacobian variety of a genus $g$  algebraic curve, one can embed  $\textrm{Kum}(T)$ into the projective space $\mathbb{P}^{2^g-1}$ to obtain a self-dual subvariety of degree $2^{g-1}g!$ with $2^{2g}$ singular points  and $2^{2g}$ trope hyperplanes, each   intersects it with multiplicity $2$ along a subvariety isomorphic to the Kummer variety of a $(g-1)$-dimensional principally polarized abelian variety.

 There is also a \emph{generalized Kummer variety} of dimension $2r$ introduced by Arnaud Beauville \cite{Beauville}.  \marginpar{Arnaud Beauville:1983} It is a nonsingular compact holomorphic symplectic manifold birationally isomorphic to the kernel of the addition map  $A^{(r+1)}\to A$, where $A$ is an abelian variety of dimension $2$ and $A^{(r+1)} = A^r/\mathfrak{S}_{r+1}$ is its symmetric product. It is one of a few known examples of complete families  of compact holomorphic symplectic manifolds of dimension $2r$. 

The first book entirely devoted to Kummer surfaces that combines geometric, algebraic and transcendental approaches to their study was published by Ronald Hudson  in 1905 \cite{Hudson}. \marginpar{Ronald Hudson:1905}

It seems that Kummer's interests to $16$-nodal quartic  surfaces arose from his pioneering study of $2$-dimensional families (congruences) of lines in $\mathbb{P}^3$ \cite{Kummer3}.  In \cite{Kummer4} (see also \cite{Kummer2}) he gives a classification of quadratic line congruences (i.e. congruences of lines such that through a general point $x\in \mathbb{P}^3$ passes exactly two lines from the congruence).  \marginpar{Ernst\,Kummer:1965}
 Let $n$ be the class of the congruence, i.e. the number of lines of the congruence that lie in a general plane. Kummer had shown that $2\le n\le 7$ and when $n= 2$ all the lines are tangent to a Kummer quartic surface at two points. There are 6 congruences like that  whose lines are tangent to the same Kummer surface.

In 1870, Felix Klein \marginpar{Felix\,Klein:1870} in his dissertation develops a beautiful relationship between the Kummer surfaces and quadratic line complexes \cite{Klein}. A quadratic line complex is the intersection $X = G\cap Q$ of the Grassmann quadric $G = \textrm{Gr}(2,4) $  in the Pl\"ucker space $\mathbb{P}^5$ with another quadric hypersurface. For any point $x\in \mathbb{P}^3$ the set of lines passing through $x$ is a plane $\sigma_x$ contained in $G$. Its intersection with $X$ is a conic.  The locus of points $x$ such that this conic becomes  reducible  is the \emph{singular surface} of the complex, and Klein had shown that it is a Kummer surface if the intersection $G\cap Q$ is transversal. Its $16$ nodes correspond to points where $\sigma_x\cap X$ is a double line. In fact, Klein shows that the set of lines in $X$ is parameterized by $\textrm{Jac}(C)$, where $C$ is the Riemann surface of genus 2, the double cover of $\mathbb{P}^1$ realized as the pencil of quadrics spanned by $G$ and $Q$ ramified along 6 points corresponding to $6$ singular quadrics in the pencil. The Kummer surface is isomorphic to $\textrm{Kum}(\textrm{Jac}(C))$. By a simultaneous diagonalization of quadrics $G$ and $Q$ he shows that the Kummer surface admits a birational nonsingular model as a complete intersection of $3$ quadrics
$$\sum_{i=1}^6 x_i^2 = \sum_{i=1}^6 a_ix_i^2 = \sum_{i=1}^6 a_i^2x_i^2 = 0,$$
where 
$y^2 = (x-a_1)\cdots (x-a_6)$ is the genus $2$ Riemann surface $C$ from above. Instead of $16$ nodes and $16$ trope-conics we now have two sets of $16$ skew lines that form a configuration $(16_6)$. Klein also shows that the Fresnel's wave surface is characterized by the condition that $C$ is bielliptic, i.e. it admits an involution with quotient an elliptic curve. By degenerating the quadratic complex to a \emph{tetrahedral quadratic complex} one obtains the  wave quartic surface. Jessop's book \cite{Jessop} gives an exposition of the works of Kummer and Klein on the relationship between line geometry and Kummer surfaces.

Klein's equations of a Kummer surface exhibit obvious symmetry defined by changing signs of the unknowns. They form an elementary abelian $2$-group $2^5$, the direct sum of $5$ copies of the cyclic group of order $2$. The quartic model also admits 16 involutions  $t_i$ defined by the projection from the nodes $p_i$ (take a general point $x$, join it with the node $p_i$, and then define $t_i(x)$ to be the residual intersection point).  In 1886 \marginpar{Felix Klein:1886} Klein asked whether the group of birational automorphisms $\textrm{Bir}(X)$ of a (general) Kummer surface is generated by the group $2^5$ and the projection involutions \cite{Klein2}. 

In 1901 John Hutchinson, \marginpar{John\,Hutchinson:1901} using the theory of theta functions, showed that a choice of one of 60  G\"opel tetrad of nodes  (means no trope-conics in its faces) leads to an equation of the Kummer surface of the form
$$q(x_1x_2+x_3x_4,x_1x_3+x_2x_4,x_1x_4+x_2x_3)+cx_1x_2x_3x_4 = 0,$$
where $q$ is a quadratic form in 3 variables \cite{Hutchinson2}. He observed that the transformation 
$x_i \mapsto 1/x_i$ leaves  this equation invariant, and hence defines a birational involution of the Kummer surface. He proves  that these 60 transformations generate  an infinite discontinuous group. However  he did not address the question whether adding these new transformations to Klein transformations would generate the whole group $\textrm{Bir}(X)$. 

In 1850 Thomas Weddle, \marginpar{Thomas\,Weddle:1850} correcting a mistake of M. Chasles, noticed that the locus of singular points of quadric surfaces passing through a fixed set of 6 general points in $\mathbb{P}^3$ is a quartic surface with singular points at the six points \cite{Weddle}. It contains the unique twisted cubic curve through the six nodes (it seems that Chasles asserted that there is nothing else). In 1861  Cayley \marginpar{Arthur Cayley:1861} gave an explicit equation of the Weddle surface \cite{Cayley2}. The simplest equation was later given by Hutchinson \cite{Hutchinson1}:
$$\det\begin{pmatrix}xyz&w&a&\alpha\\
yzw&x&b&\beta\\
xzw&y&c&\gamma\\
xyw&z&d&\delta\end{pmatrix} = 0
$$
Here the six points are the reference points $(1:0:0:0),(0:1:0:0),(0:0:1:0),(0:0:0:1)$ and two points with coordinates given in the last two columns of the matrix. If one inverts the coordinates, the equation does not change. In 1889 Schottky \marginpar{Friedrich Schottky:1889} proved that a Weddle surface is birationally isomorphic to a Kummer surface \cite{Schottky}. This implies that the groups of birational automorphisms of a Weddle surface and a Kummer surface are isomorphic. For example, it is immediate to see that the same inversion transformation used by Hutchinson leaves the Weddle surface invariant. In 1911 Virgil Snyder \marginpar{Virgil\,Snyder:1911} gave many geometric constructions of  involutions of the Weddle surface that give corresponding involutions of the Kummer surface \cite{Snyder}. He somehow does not mention \emph{Hutchinson involution}.\footnote{John Hutchinson, Virgil Snyder, Francis Sharpe and Clide Craig were on the faculty of the Cornell Mathematics Department.}   In a paper of 1914, Francis Sharpe and Clide Craig \marginpar{Francis\,Sharpe:1914} pioneered the new approach to study birational automorphisms of algebraic surfaces based on Francesco Severi's Theory of a Basis. \marginpar{Clide\,Craig:1914} It consists of representing an automorphism of an algebraic surface as a transformation of the group of algebraic cycles on the surface. This gave an easy proof of Hutchinson's result that 60 G\"opel-Hutchinson involutions generate an infinite group \cite{Sharpe}.

The Hutchinson involutions act freely on a nonsingular model of the Kummer surface and the quotient by these involutions are \emph{Enriques surfaces}. \marginpar{Shigeru\,Mukai:2013}\marginpar{Hasinori\,Ohashi:2013} In their paper of 2013 
under the title `Enriques surfaces of Hutchinson-G\"opel type and Mathieu automorphisms' Shigeru Mukai and Hasinori Ohashi study the automorphism group of these Enriques surfaces \cite{Mukai}. 

John Hutchinson also discovered the following amazing fact.\marginpar{John\,Hutchinson:1899} The linear system of quadric surfaces through a \emph{Weber hexad} of nodes on a Kummer surface defines a birational isomorphism to the \emph{Hessian surface} of a nonsingular cubic surface $F_3(x,y,z,w) = 0$ (the Hessian surface defined by the hessian matrix of $F_3(x,y,z,w)$) \cite{Hutchinson1}. The Hessian surface can be also defined as the locus of singular \emph{polar quadrics} of the surface (quadrics given by linear combinations of partial derivatives of $F_3$). The rational map that assigns to a singular quadric its singular point defines a fixed-point-free involution on a nonsingular model of the Hessian surface with quotient isomorphic to an Enriques surface. There are 120 Weber hexads, each defines a birational involution of a Kummer surface, called the \emph{Hutchinson-Weber involution}.  
The automorphism groups as well as complex dynamics of the Hessian quartic surfaces are the subjects of study of several recent papers.

It is not surprising that none of the classical geometers could decide whether a given finite set of transformations generates the group $\textrm{Bir}(X)$.  It had to wait until the end of the century for  new technical tools to arrive.

A  nonsingular birational model of a Kummer surface is an example of a K3 surface. By definition, a K3 surface $Y$ is a compact analytic simply-connected surface with trivial first Chern class $c_1(Y)$. Its second Betti number is equal to $22$. By a theorem of John Milnor,  its homotopy type   is uniquely determined by the quadratic form expressing the cup-product on $H^2(Y,\mathbb{Z}) \cong \mathbb{Z}^{22}$.  It is a unique unimodular even quadratic lattice of signature $(3,19)$ isomorphic to $E_8(-1)^{\oplus 3}\oplus U$, where $E_8(-1)$ is an even unimodular negative definite lattice of rank $8$ and $U$ is the integral hyperbolic plane. All K3 surfaces are diffeomorphic and realize the homotopy type defined by this  lattice.   Any complex K3 surface  admits a K\"ahler metric and, by Yau's theorem it also admits a \emph{Ricci-flat metric}. It is not known how to write it explicitly (a question of great importance for physicists). The Kummer surface is flat in this metric outside the singular points, but still one does not know how to extend this flat metric to a Ricci-flat metric on $Y$. 

The quadratic lattice of cohomology $H^2(Y,\mathbb{Z})$ of an  algebraic K3 surface contains   a sublattice 
$S_Y = H^2(Y,\mathbb{Z})_{alg} \cong \mathbb{Z}^{\rho}$ of algebraic $2$-cycles. It contains  $c_1(L)$, where $L$ is an ample line bundle. The signature of $S_Y$ is equal to $(1,\rho-1)$ and, in general, it is not unimodular. The group $\textrm{Bir}(X)$ is isomorphic to the group $\textrm{Aut}(Y)$ of biregular automorphisms of $Y$ and it admits a natural representation $\rho:\textrm{Aut}(Y)\to \text{O}(H^2(Y,\mathbb{Z}))$ in the orthogonal group of $H^2(Y,\mathbb{Z})$ that leaves $S_Y$ invariant. 

The fundamental \emph{Global Torelli Theorem} of I. Pyatetsky-Shapiro and I. Shafarevich \cite{PSSh}\marginpar{Ilya\,Pyatetsky-Shapiro:1977} \marginpar{Igor\,Shafarevich:1977}
allows one to describe the image of  homomorphism $\rho$ as the set of all isometries of $H^2(Y,\mathbb{Z})$ that leave (after complexification) $H^{2,0}(Y,\mathbb{C})$ invariant and also leave invariant the semi-group of cohomology classes of holomorphic curves  on $Y$. The group of automorphisms $\textrm{Aut}(Y)$ acts naturally on the lattice $S_Y$ and on the \emph{hyperbolic space} $\mathbb{H}^{\rho-1}$ associated with the linear space $S_Y\otimes \mathbb{R}$ of signature $(1,\rho-1)$. In this way it is realized as a discrete group of motions of a hyperbolic space. By using an isometric embedding of $S_Y$ into the unimodular even lattice $II_{1,25}$ of signature $(1,25)$ isomorphic to the orthogonal sum of the \emph{Leech lattice} and the hyperbolic plane $U$,  Richard Borcherds \marginpar{Richard Borcherds:1997} introduced a method that in some cases allows one to compute the automorphism group of a K3 surfaces using the isometries of $II_{1,25}$ defined by the reflections into \emph{Leech roots} \cite{Borcherds}. In 1998, based on Borcherds' ideas,  Shigeyuki Kondo \marginpar{Shigeyuki\,Kondo:1998} proved  that the group of birational automorphisms of a general Kummer surface (i.e. the Jacobian surface of a general curve of genus 2) is generated by the group $2^5$, 16 projection involution, 60 Hutchinson-G\"opel involutions and 120 Hutchinson-Weber involutions \cite{Kondo}.\footnote{In fact,  Kondo includes in the set of generators  some new automorphisms of infinite order discovered by  JongHae Keum \cite{Keum} instead of Hutchinson-Weber involutions.  It was later observed by Ohashi that one can use Hutchinson-Weber involutions instead of Keum's automorphisms \cite{Ohashi}} 
 The group of birational automorphisms of an arbitrary Kummer surface is still unknown.

 We refer to some modern expositions of the theory of Kummer surfaces to \cite{Lange},  \cite{Dolgachev} and \cite{GH}. We apologize for the brevity of our discussion that omits many important contributions to the study of Kummer surfaces in the past and in the present. Note that a search of the  data base of Math.Sci.Net under Title:Kummer surface  gives 130 items (70 of them in this century).

\end{document}